\documentclass[11pt,a4paper,leqno]{amsart}
\usepackage{amssymb,graphicx}
\usepackage{hyperref}
\makeatletter\renewcommand{\@secnumfont}{\relax}\makeatother
\theoremstyle{plain}

\newtheorem{thm}{Theorem}[section]

\newtheorem{conj}[thm]{Conjecture}
\theoremstyle{definition}

\theoremstyle{remark}

\newcommand{\proofof}[1]{\end{#1}\begin{proof}} 
\numberwithin{equation}{section}

\newcommand{\thismonth}{\ifcase\month\or
  January\or February\or March\or April\or May\or June\or
  July\or August\or September\or October\or November\or December\fi
  \space\number\year}
\newcommand{\sideremark}[1]{\ifvmode\leavevmode\fi\vadjust
 {\vbox to0pt{\vss\hbox to0pt{\hskip\hsize\hskip1em\vbox{\hsize3cm\tiny
  \raggedright\pretolerance10000\noindent
  #1\hfill}\hss}\vbox to8pt{\vfil}\vss}}}
\DeclareSymbolFont{script}{U}{eus}{m}{n}
\DeclareSymbolFontAlphabet{\mathscr}{script}
\DeclareMathSymbol{\EuWedge}{0}{script}{"5E}
\DeclareMathAlphabet{\mathrmsl}{OT1}{cmr}{m}{sl}
\newcommand{\rssymb}[2]{\newcommand{#1}{{\mathrmsl{#2}}}}
\newcommand{\calsymb}[2]{\newcommand{#1}{{\mathcal{#2}}}}
\newcommand{\bbsymb}[2]{\newcommand{#1}{{\mathbb{#2}}}}
\newcommand{\lieoper}[2]{\newcommand{#1}{\mathop{\mathfrak{#2}}}}
\newcommand{\oper}[3][n]{\newcommand{#2}{\mathop{\mathrm{#3}}\ifx
  n#1\nolimits\else\limits\fi}}
\newcommand{\rsoper}[3][n]{\newcommand{#2}{\mathop{\mathrmsl{#3}}\ifx
  n#1\nolimits\else\limits\fi}}
\bbsymb\C{C}\bbsymb\F{F}\bbsymb\HQ{H}\bbsymb\N{N}\bbsymb\Q{Q}
\bbsymb\R{R}\bbsymb\U{U}\bbsymb\V{V}\bbsymb\W{W}\bbsymb\Z{Z}
\calsymb\cA{A}\calsymb\cB{B}\calsymb\cC{C}\calsymb\cD{D}\calsymb\cE{E}
\calsymb\cF{F}\calsymb\cG{G}\calsymb\cH{H}\calsymb\cI{I}\calsymb\cJ{J}
\calsymb\cK{K}\calsymb\cL{L}\calsymb\cM{M}\calsymb\cN{N}\calsymb\cO{O}
\calsymb\cP{P}\calsymb\cQ{Q}\calsymb\cR{R}\calsymb\cS{S}\calsymb\cT{T}
\calsymb\cU{U}\calsymb\cV{V}\calsymb\cW{W}\calsymb\cX{X}\calsymb\cY{Y}
\calsymb\cZ{Z}
\renewcommand{\geq}{\geqslant}\renewcommand{\leq}{\leqslant}
\oper\End{End}                       
\oper\Sym{Sym}\oper\Skew{Skew}
\oper\Cl{Cl}
\oper\Aut{Aut}                       
\oper\GL{GL}\oper\SL{SL}
\oper\CO{CO}\oper\Orth{O}\oper\SO{SO}\oper\Pin{Pin}\oper\Spin{Spin}
\oper\Symp{Sp}\oper\Un{U\strut}\oper\SU{SU}
\rsoper\Diff{Diff}\rsoper\SDiff{SDiff}
\lieoper\der{der}                    
\lieoper\gl{gl}\lieoper\sgl{sl}
\lieoper\co{co}\lieoper\orth{o}\lieoper\so{so}\lieoper\spin{spin}
\lieoper\symp{sp}\lieoper\un{u}\lieoper\su{su}
\rsoper\Vect{Vect}
\rsoper\kernel{ker}\rsoper\image{im} 
\rsoper\alt{alt}\rsoper\sym{sym}     
\oper\Ad{Ad}\rsoper\ad{ad}           
\newcommand{\restr}[1]{|_{\lower1pt\hbox{${}_{#1}$}}}
\newcommand{\norm}[2][]{|\mkern-2mu|#2|\mkern-2mu|_{\lower1pt\hbox{${}_{#1}$}}}

\newcommand{\inner}{\mathbin{\raise1pt\hbox{$\lrcorner$}}}
\newcommand{\capinner}{\mathbin{\lower3pt\hbox{$\urcorner$}}}

\def\st.{\mathrel{|}}
\def\kahl/{k\"ahler}
\oper\trace{tr} 
\rsoper\divg{div} 
\rssymb\iden{id}  
\rssymb\vol{vol}  
\rssymb\ev{ev}    
\rsoper\scal{scal}

\DeclareMathOperator{\Obs}{Obs}

\newcommand{\sd}{{\raise1pt\hbox{$\scriptscriptstyle +$}}}
\newcommand{\asd}{{\raise1pt\hbox{$\scriptscriptstyle -$}}}
\newcommand{\sdasd}{{\raise1pt\hbox{$\scriptscriptstyle\pm$}}}
\newcommand{\asdsd}{{\raise1pt\hbox{$\scriptscriptstyle\mp$}}}

\def\Dbar{\leavevmode\lower.6ex\hbox to 0pt{\hskip-.23ex \accent"16\hss}D}

\begin{document}
\title{Configuration spaces of points, symmetric groups and polynomials of several variables}
\author{Joseph Malkoun}
\address{Department of Mathematics and Statiscs\\
Notre Dame University-Louaize\\
Lebanon}
\email{joseph.malkoun@ndu.edu.lb}
\date{\today}
\begin{abstract}
Denoting by $C_n(X)$ the configuration space of $n$ distinct points in $X$, with $X$ being either Euclidean $3$-space $\mathbb{E}^3$ or 
hyperbolic $3$-space $\mathbb{H}^3$ or $\mathbb{C}P^1$ , by 
$\mathscr{P}_{k,d}$ the vector space of homogeneous complex polynomials in the variables $z_0, \ldots, z_k$ of degree $d$, and by 
$\Obs^n_d$ the set of all $d$-subsets of $\{1,\ldots,n\}$, the symmetric group $\Sigma_n$ acts on $C_n(\mathbb{R}^3)$ by 
permuting the $n$ points and also acts in a natural way on $\Obs^n_d$. With $n = k+d$, the space $\mathscr{P}_{k,d}$ has dimension $\binom{n}{d}$, 
which is also the number of elements in $\Obs^n_d$. It is thus natural to ask the following question. Is there a family of 
continuous maps $f_I: C_n(X) \to \mathbb{P}\mathscr{P}_{k,d}$, for $I \in \Obs^n_d$ (here $\mathbb{P}$ is complex projectivization), 
which satisfies $f_I(\sigma.\mathbf{x}) = f_{\sigma.I}(\mathbf{x})$, for all $\sigma \in \Sigma_n$ and all $\mathbf{x} \in C_n(X)$, 
and such that, for each $\mathbf{x} \in C_n(X)$, the polynomials $f_I(\mathbf{x})$, for $I\in \Obs^n_d$, each 
defined up to a scalar factor, are linearly independent over $\mathbb{C}$?  We 
provide two closely related smooth candidates for such maps for each of the two cases, Euclidean and hyperbolic, which would be solutions to the 
above problem provided a linear independence conjecture holds. 
Our maps are natural extensions of the Atiyah-Sutcliffe maps. Moreover, by taking in the hyperbolic versions the limiting case of $n$ distinct points on the 
sphere at infinity, thought of as the Riemann sphere, we get two constructions of actual solutions of the above problem for $X = \mathbb{C}P^1$, as we 
prove linear independence for these last two constructions. These last two constructions are classical in character, and can be viewed as higher dimensional versions of 
Lagrange polynomial interpolation. They appear to be new, though we suspect they are not. 
\end{abstract}

\maketitle

\section{Introduction}

The origin of the idea is rooted in the Atiyah-Sutcliffe paper \cite{Atiyah-Sutcliffe2002}. We briefly summarize some key points of that paper. We denote by $C_n(\mathbb{R}^3)$ the manifold consisting of all $n$-tuples of distinct points in $\mathbb{R}^3$. The Berry-Robbins problem asks whether there exists, for any $n \geq 2$, a continuous map
\[ f_n: C_n(\mathbb{R}^3) \to U(n)/T \]
where $T$ is the diagonal maximal torus of $U(n)$, which is equivariant under the action of the symmetric group $\Sigma_n$. The latter permutes the $n$ points $x_1, \ldots x_n$ of a configuration $\mathbf{x} = (x_1, \ldots, x_n) \in C_n(\mathbb{R}^3)$, and permutes the $n$ columns of an element $g \in U(n)/T$. The Berry-Robbins problem can be relaxed a little. It is enough to find, for each $n \geq 2$, a continuous map
\[ F_n: C_n(\mathbb{R}^3) \to GL(n,\mathbb{C})/(\mathbb{C}^*)^n \]
which is equivariant under the action of the symmetric group $\Sigma_n$.

While it was answered positively in \cite{Atiyah2000}, the maps constructed there are not smooth. On the other hand, in \cite{Atiyah-Bielawski2002}, the authors generalize the Berry-Robbins problem for any compact Lie group, and solve it in that general setting. The maps they obtain are not explicit though, since they rely on an analysis of Nahm's equations.

On the other hand, in the papers \cite{Atiyah2000} and \cite{Atiyah2001}, Sir Michael Atiyah constructs smooth candidates for solutions to the Berry-Robbins problem, which would be genuine solutions provided a linear independence conjecture holds.

The idea is as follows. Given a configuration $\mathbf{x} = (x_1, \ldots, x_n) \in C_n(\mathbb{R}^3)$, 
from each point $x_i$, and for each $j \neq i$, form the points
\[ t_{ij} = \frac{x_j-x_i}{\norm{x_j-x_i}} \in S^2 \]
Then use stereographic projection to identify $S^2$ with the Riemann sphere $\hat{\mathbb{C}} = 
\mathbb{C} \cup \{ \infty \}$. Then, for each $i$, $1 \leq i \leq n$, form the polynomial $p_i$ having the points $t_{ij}$ for $j \neq i$, as roots. Such a polynomial $p_i$ is uniquely determined up to a non-zero scalar factor and, as a complex polynomial, has degree at most $n-1$. Sir Michael Atiyah conjectures the following.
\begin{conj} \label{AC} For any configuration $\mathbf{x}$, the corresponding polynomials $p_1$, \ldots, $p_n$ are linearly independent over $\mathbb{C}$.\end{conj} 
This conjecture was proved for $n=3$ by Sir Michael Atiyah in \cite{Atiyah2000} and \cite{Atiyah2001}, and for $n=4$ by Eastwood and Norbury in \cite{Eastwood-Norbury2001}, and for some special 
configurations by {\Dbar}okovi{\'c} in \cite{Dokovic2002a} and \cite{Dokovic2002b}.

Moreover, in \cite{Atiyah-Sutcliffe2002}, the authors construct a normalized determinant function 
\[ D: C_n(\mathbb{R}^3) \to \mathbb{C} \]
whose non-vanishing is equivalent to the previous conjecture \ref{AC}. Moreover, they formulate what they call conjectures $2$ and $3$ (with conjecture $3$ implying conjecture $2$). We shall not discuss conjecture $3$, but only conjecture $2$, very briefly. 

\begin{conj}[Atiyah-Sutcliffe conjecture 2] \label{ASC2} $|D(\mathbf{x})| \geq 1$ for any configuration $\mathbf{x} \in 
C_n(\mathbb{R}^3)$. \end{conj}
This conjecture was proved in \cite{Atiyah-Sutcliffe2002} for $n=3$, and was proved for $n=4$ by Bou 
Khuzam and Johnson in \cite{BKJ2014} and using a different method by Svrtan in \cite{Svrtan2014} (in fact, also the stronger conjecture 3 is proved in \cite{BKJ2014} and \cite{Svrtan2014}).

In this paper, we generalize the Atiyah-Sutcliffe maps, as well as Conjectures \ref{AC} and 
\ref{ASC2}. Our basic idea can be summarized in the following way. Instead of partitioning the $n$ distinct points $x_1$, \ldots, $x_n$ into a single observer and $n-1$ stars, we partition them 
instead into $d$ observers and $n-d$ stars. It turns out that this immediately leads us to work with complex homogeneous polynomials of degree $d$ in $n-d+1$ complex variables, and the Atiyah-Sutcliffe normalized determinant function generalizes to our setting.

In the following, $X$ will denote either Euclidean $3$-space $\mathbb{E}^3$, or hyperbolic $3$-space $\mathbb{H}^3$, or $\mathbb{C}\mathbb{P}^1$. Denoting by $\mathscr{P}_{k,d}$ the vector space of homogeneous complex polynomials in the variables $z_0, \ldots, z_k$ of degree $d$, and by 
$\Obs^n_d$ the set of all $d$-subsets of $\{1,\ldots,n\}$, the symmetric group $\Sigma_n$ acts on $C_n(\mathbb{R}^3)$ by 
permuting the $n$ points and also acts in a natural way on $\Obs^n_d$. With $n = k+d$, the space $\mathscr{P}_{k,d}$ has dimension $\binom{n}{d}$, 
which is also the number of elements in $\Obs^n_d$. It is thus natural to ask the following question. Is there a family of 
continuous maps $f_I: C_n(X) \to \mathbb{P}\mathscr{P}_{k,d}$, for $I \in \Obs^n_d$ (here $\mathbb{P}$ is complex projectivization), 
which satisfies $f_I(\sigma.\mathbf{x}) = f_{\sigma.I}(\mathbf{x})$, for all $\sigma \in \Sigma_n$ and all $\mathbf{x} \in C_n(\mathbb{R}^3)$, 
and such that, for each $\mathbf{x} \in C_n(\mathbb{R}^3)$, the polynomials $f_I(\mathbf{x})$, for $I\in \Obs^n_d$, each 
defined up to a scalar factor, are linearly independent over $\mathbb{C}$?

We present two closely related constructions of maps, which are smooth candidates for solutions of the previous problem for $X = \mathbb{E}^3$ and $X=\mathbb{H}^3$, and are genuine solutions 
provided a linear independence conjecture is true, completely similar to the Atiyah-Sutcliffe case (which our construction generalizes). We also generalize the construction of the Atiyah-Sutcliffe 
normalized determinant function, and generalize Conjecture \ref{ASC2}. In the case where $X=\mathbb{C}\mathbb{P}^1$, we obtain two elegant solutions, for which we are able to prove linear independence. These latter two constructions arise from the two hyperbolic constructions by taking the limiting case as the configuration tends to a configuration of $n$ distinct points on the sphere 
at infinity, which we think of as the Riemann sphere or, in other words, as $\mathbb{C}\mathbb{P}^1$. These two constructions, at least as presented here, appear to be new, 
though we suspect that they are not, being very classical in nature, and natural generalizations of the Lagrange polynomial interpolation problem to polynomials of several variables.

\section{Observer-based maps (first construction)}

Fix an integer $n \geq 2$, and an integer $d$, with $1 \leq d \leq n-1$. Let $z_0,\ldots,z_k$ 
be complex variables, where $k=n-d$. For convenience, we use the multi-index notation. A multi-index for the 
variables $z_0,\ldots,z_k$ is a $k+1$-tuple $M=(m_0,\ldots,m_k) \in \mathbb{N}^{k+1}$, and we define $\mathbf{z}^M = z_0^{m_0}\ldots z_k^{m_k}$. 
The length $|M|$ of $M$ is defined by $|M| = m_0 + \cdots + m_k$. We present here the Euclidean version of the observer-based maps, and indicate at the end of 
this section how to modify the construction to get a hyperbolic version, as well as a Riemann sphere version. We introduce the configuration 
spaces $C_n(\mathbb{R}^3)$, and the polynomial spaces $\mathscr{P}_{k,d}$.
\begin{align*}
C_n(\mathbb{R}^3) &= \{ \mathbf{x} = (\mathbf{x}_1,\ldots,\mathbf{x}_n) \in (\mathbb{R}^3)^n; \mathbf{x}_i \neq \mathbf{x}_j, 
\text{ for $i\neq j$, $1 \leq i,j \leq n$} \} \\
\mathscr{P}_{k,d} &= \{ \sum_{|M|=d} c_M \mathbf{z}^M; c_M \in \mathbb{C} \}
\end{align*} 
We denote by $\Obs^n_d$ to be the set of all $d$-subsets of $\{1,\ldots,n\}$. Fix $I \in \Obs^n_d$. 
Write $I = \{ i_1, \ldots, i_d \}$, where $1 \leq i_1 < \ldots < i_d \leq n$, without loss of generality. 
Given a configuration $\mathbf{x}=(\mathbf{x}_1, \ldots, \mathbf{x}_n) \in C_n(\mathbb{R}^3)$, 
we think of a choice of $I \in \Obs^n_d$ as a choice of $d$ ``observers'' among the $n$ distinct points $\mathbf{x}_1,\dots,\mathbf{x}_n$ 
(namely, the $d$ points $\mathbf{x}_{i_1},\ldots,\mathbf{x}_{i_d}$). This explains the 
notation $\Obs^n_d$, since it is the set of all possible choices of $d$ observers among the $n$ points. Given a 
choice of $I \in \Obs^n_d$, we consider $J = I^c$, the complement of $I$ in the 
set $\{1,\dots,n \}$. Then $|J| = n-d = k$. We think of $J$ as corresponding to a set of $k$ ``stars'' among the $n$ 
points $\mathbf{x}_1,\ldots,\mathbf{x}_n$, namely the points $\mathbf{x}_{j_1},\ldots,\mathbf{x}_{j_k}$, where 
$J = \{ j_1,\ldots, j_{k} \}$, and we assume without loss of generality that $1 \leq j_1 < \ldots < j_k \leq n$. Summarizing, a choice 
of $I \in \Obs^n_d$ corresponds to partitioning the $n$ points $\mathbf{x}_1,\ldots,\mathbf{x}_n$ into $d$ observers 
$\mathbf{x}_{i_1},\dots,\mathbf{x}_{i_d}$ and $k$ stars $\mathbf{x}_{j_1},\dots, \mathbf{x}_{j_{k}}$, with $J = I^c$.

Consider $\mathbb{C}^2$ with its complex coordinates $u$ and $v$. 
Then $S^3 = \{ (u,v)\in \mathbb{C}^2; |u|^2 + |v|^2 = 1 \}$, and we define the Hopf map $h: S^3 \to S^2$ by
\[ h(u,v) = (2\bar{u}v, |v|^2-|u|^2) \]
where $S^2$ is the unit sphere in $\mathbb{R}^3$, itself thought of as $\mathbb{C} \times \mathbb{R}$. For each 
pair $(i,j)$, $1 \leq i,j \leq n$ and $i \neq j$, we choose a Hopf lift $(u_{ij},v_{ij}) \in S^3 \subset \mathbb{C}^2$ of
\begin{equation} \label{tij} t_{ij} = \frac{\mathbf{x}_j-\mathbf{x}_i}{\norm{\mathbf{x}_j-\mathbf{x}_i}} \in S^2 \end{equation}
where $\norm{\,}$ denotes the Euclidean norm. We assume that once a choice of Hopf lift for the pair $(i,j)$, with $i<j$, 
is made, the Hopf lift corresponding to $(j,i)$ is
\[ (u_{ji},v_{ji}) = (-\bar{v}_{ij}, \bar{u}_{ij}) \]
This is valid since the point $(\mathbf{x}_i-\mathbf{x}_j)/\norm{\mathbf{x}_i-\mathbf{x}_j}\in S^2$ corresponding to 
$(j,i)$ is the antipodal of the point $(\mathbf{x}_j-\mathbf{x}_i)/\norm{\mathbf{x}_j-\mathbf{x}_i}\in S^2$  corresponding to $(i,j)$. 
Once the choices of Hopf lifts are made, we can associate to each pair $(i,j)$, for $i \neq j$, a homogeneous complex 
polynomial $L_{ij}$ depending on $u$ and $v$ in the following way.
\begin{align*} L_{ij}(u,v) &= \left| \begin{array}{cc} u_{ij} & u \\
v_{ij} & v \end{array} \right| \\
&= u_{ij}v-v_{ij}u
\end{align*}
Given a choice of $I \in \Obs^n_d$, and for each $i \in I$, form the homogeneous polynomial
\[ q_i(u,v) = \prod_{j \in I^c} L_{ij}(u,v) \in \mathscr{P}_{1,k} \]
of degree $k$ in $u$ and $v$. Consider two arbitrary elements $q$ and $\tilde{q}$ of $\mathscr{P}_{1,k}$, which can be written as
\begin{align*} 
q(u,v) &= c_0 u^k + c_1 u^{k-1} v + \cdots + c_{k-1} u v^{k-1} + c_k v^k \\
\tilde{q}(u,v) &= d_0 u^k + d_1 u^{k-1} v + \cdots + d_{k-1} u v^{k-1} + d_k v^k
\end{align*}
Define the following non-degerate complex bilinear form $(-,-)$ on $\mathscr{P}_{1,k}$ by
\[ (q,\tilde{q}) = \sum_{i=0}^k (-1)^i \frac{c_i d_{k-i}}{\binom{k}{i}} \]
Given an element $q \in \mathscr{P}_{1,k}$ denote by $\Lambda_q \in \mathscr{P}^*_{1,k}$ the following
\[ \Lambda_q(-) = (q,-) \]

We can now define (given $I \in \Obs^n_d$)
\[ p_I = \bigodot_{i \in I}  \Lambda_{q_i} \]
Here $\bigodot$ denotes the symmetric tensor product. Given a choice of Hopf lifts, for each $I \in \Obs^n_d$, $p_I \in \mathscr{P}_{k,d}$. Thus, 
we have constructed, for each $I \in \Obs^n_d$, a smooth map $f_I: C_n(\mathbb{R}^3) \to P(\mathscr{P}_{k,d})$, mapping a configuration $\mathbf{x}$ to 
$[p_I] \in P(\mathscr{P}_{k,d})$, where $[ -]$ denotes a non-zero polynomial up to (non-zero) constant scaling. Moreover, for every $\sigma \in \Sigma_n$, 
where $\Sigma_n$ is the permutation group on $n$ elements, we have the following equivariance property
\[ f_I(\sigma.\mathbf{x}) = f_{\sigma.I}(\mathbf{x}) \]
We now make the following conjecture:
\begin{conj}[Linear Independence Conjecture for the observer-based maps] \label{conj1} For each configuration $\mathbf{x} \in C_n(\mathbb{R}^3)$, the corresponding polynomials $p_I$, as $I$ varies in $\Obs^n_d$, are linearly independent over $\mathbb{C}$. \end{conj}

We now indicate how to modify the previous construction for $X = \mathbb{H}^3$. Instead of the formula \ref{tij}, the point $t_{ij}$ is in this case defined to be the limiting point on the sphere 
at infinity, thought of as the Riemann sphere, in the Poincare ball model of $\mathbb{H}^3$, of the hyperbolic ray emanating from point $x_i$ and passing through point $x_j$. The construction of the 
polynomials $p_I$ then proceeds without further modifications.

\section{Star-based maps (second construction)}

This is a closely related construction. We now assume we have $k$ observers and $d$ stars, thus interchanging the roles of $k$ and $d$. Our discussion will be for either $X = \mathbb{E}^3$ or 
$X = \mathbb{H}^3$. We first proceed just as in our first construction of observer-based maps, by choosing Hopf lifts for the $t_{ij}$. Bear in mind that the $t_{ij}$ are given by formula 
\ref{tij} in the Euclidean case, but in the hyperbolic case, $t_{ij}$ is the point at infinity towards which tends the hyperbolic ray from $x_i$ passing through $x_j$. Given $I \in \Obs^n_k$, we let 
$J = I^c$ (the complement of $I$), and now for each $j \in J$, we form
\[ q_j(u,v) = \prod_{i \in I} L_{ij}(u,v) \]
Thus now, and unlike in our first construction, a star is fixed, and we take the product over all observers, which explain the terminology ``star-based''. 
Thus each $q_j$ is an element of $\mathscr{P}_{1,k}$. We then form
\[ p_I = \bigodot_{j \in J} \Lambda_{q_j} \]
We also make a similar conjecture for the star-based maps.
\begin{conj}[Linear independence conjecture for the star-based maps] \label{conj2} For each configuration $\mathbf{x} \in C_n(\mathbb{R}^3)$, the 
corresponding polynomials $p_I$, as $I$ varies in $\Obs^n_k$, are linearly independent over $\mathbb{C}$. \end{conj}

\section{Normalized Determinant functions}
We now wish to define a normalized determinant function, similar to the Atiyah-Sutcliffe normalized determinant function in \cite{Atiyah-Sutcliffe2002} (and in fact generalizes the Atiyah-Sutcliffe normalized determinant, which corresponds to $d=1$ in our observer-based maps). We explain how to define such normalized determinant functions for the observer-based construction only, the star-based construction being similar. 

Given $I \in \Obs^n_d$, where $I = \{ i_1, \ldots, i_d \}$, and $1 \leq i_1 < \ldots < i_d \leq n$, we consider $J = I^c$, where $J= \{ j_1, \dots, j_k \}$ with $1 \leq j_1 < \ldots < j_k \leq n$, and associate to $I$ the following monomial of degree $d$:
\[ q_I(z_0,\cdots, z_k) = z_0^{j_1-1} z_1^{j_2-j_1-1}\cdots z_{j_l}^{j_{l+1}-j_l-1} \cdots z_{k-1}^{j_k-j_{k-1}-1} z_k^{n-j_k} \]
It is clear that the $q_I$, as $I$ varies in $\Obs^n_d$, form a basis of $V_{k,d}$. 

We now assume that $\Obs^n_d$ is endowed with the lexicographic order. This also induces an order on the basis $(q_I)$ that we have just defined. Thus given a configuration $\mathbf{x} \in C_n(X)$, where $X$ is either $\mathbb{E}^3$ or $\mathbb{H}^3$, and once the choices of Hopf lifts are made, with the only restriction that the Hopf lifts of $t_{ij}$ and $t_{ji}$ form a symplectic basis of $\mathbb{C}^2$, we can combine the polynomials $p_I$ into a single $\binom{n}{d}$ by $\binom{n}{d}$ matrix $M$, whose $I$'th column contains the coefficients of $p_I$ with respect to the basis $(q_{I'})$ ordered lexicographically.

We can now define the normalized determinant function $D: C_n(\mathbb{R}^3) \to \mathbb{C}$ by the simple formula:
\[ D(\mathbf{x}) = \det(M) \]
One can show that the value of $D(\mathbf{x})$ is indeed well defined, and independent of the choices of Hopf lifts made previously. Conjecture \ref{conj1} is equivalent to the non-vanishing of $D$. Similar to Atiyah-Sutcliffe's conjecture 2 in \cite{Atiyah-Sutcliffe2002}, we make the following stronger conjecture (for both the observer-based and star-based constructions).
\begin{conj} \label{conj2} $|D(\mathbf{x})| \geq 1$ for any configuration $\mathbf{x} \in C_n(X)$, where $X$ is either $\mathbb{E}^3$ or $\mathbb{H}^3$.
\end{conj}
Similar to the Atiyah-Sutcliffe setting, our function $D$ is invariant under Euclidean transformations and scaling in $\mathbb{R}^3$ (and their induced action on $C_n(\mathbb{R}^3)$), and is also invariant under the action of the symmetric group $\Sigma_n$ on $C_n(\mathbb{R}^3)$. 

\section{Versions of our two constructions for $X = \mathbb{C}\mathbb{P}^1$}

For $X = \mathbb{C}\mathbb{P}^1$, one can similarly obtain two constructions, an observer-based one, and a star-based one. Moreover, one can prove linear independence 
for these two constructions in this case. These can be obtained from the hyperbolic case, by letting the configuration tend to a configuration of $n$ distinct points on the sphere 
at infinity, thought of as the Riemann sphere. However, we prefer to present the constructions more directly.

\subsection{Observer-based construction for $X = \mathbb{C}\mathbb{P}^1$} 
Let $z_1,\ldots, z_n$ be $n$ distinct points in $\mathbb{C}\mathbb{P}^1$, and let $I \in \Obs^n_d$ be a $d$-subset of $\{1,\ldots,n\}$, and denote by $J = I^c$ its 
complement. We choose Hopf lifts of $z_1, \ldots, z_n$, which we denote by $(u_1,v_1), \ldots, (u_n,v_n)$ respectively. We introduce the linear forms $L_i(u,v)$.
\[ L_i(u,v) = u_i v - v_i u \]
We define the homogeneous polynomial $q_I \in \mathscr{P}_{1,k}$ of the two complex variables $u$ and $v$ of degree $k$ by
\[ q_I(u,v) = \prod_{j\in J} L_j(u,v) \]
We then define
\[ p_I = \Lambda_{q_I} \odot \cdots \odot \Lambda_{q_I} \in \mathscr{P}_{k,d} \]
where the number of $\Lambda_{q_I}$ factors on the right-hand side is $d$. To prove that the $p_I$ are linearly independent over $\mathbb{C}$, it suffices to exhibit 
$\binom{n}{d}$ points $g_I$ in $\bigodot^d \mathscr{P}_{1,k}$, also indexed by $I \in \Obs^n_d$, such that
\[ p_I(g_{I'}) = \left\{ \begin{array}{ll} 1, & \text{if $I = I'$} \\
                                                       0, & \text{if $I \neq I'$} \end{array} \right. \]
Such points $g_I$ are given by
\[ g_I = \bigodot_{i \in I} L_i^k \]
This is because if $q \in \mathscr{P}_{1,k}$, then
\[ (p, (u_0 v - v_0 u)^k) = p(u_0,v_0) \]
as can be directly checked.

\subsection{Star-based construction for $X = \mathbb{C}\mathbb{P}^1$}

We also start with $n$ distinct points $z_1,\ldots, z_n$ in $\mathbb{C}\mathbb{P}^1$. We now let $I \in \Obs^n_k$ be a $k$-subset of $\{1,\ldots,n\}$, and denote 
by $J = I^c$ its complement. As in the observer-based construction, we choose Hopf lifts of $z_1, \ldots, z_n$, which we denote by 
$(u_1,v_1), \ldots, (u_n,v_n)$ respectively. We introduce the linear forms $L_i(u,v)$.
\[ L_i(u,v) = u_i v - v_i u \]
For $j \in J$, we let
\[ q_j(u,v) = L_j(u,v)^k \]
and then define
\[ p_I = \bigodot_{j \in J} \Lambda_{q_j} \in \mathscr{P}_{k,d} \]
We introduce
\[ h_I(u,v) = \prod_{i \in I} L_i(u,v) \]
and then define
\[ g_I = h_I \odot \cdots \odot h_I \in \bigodot^d \mathscr{P}_{1,k} \]
It can then be checked, similar to the observer-based construction, that
\[ p_I(g_{I'}) = \left\{ \begin{array}{ll} 1, & \text{if $I = I'$} \\
                                                       0, & \text{if $I \neq I'$} \end{array} \right. \]
This establishes linear independence of the $p_I$ for the star-based construction.

\section{concluding remarks}

While the conjectures pertaining to the Euclidean and hyperbolic versions of our constructions seem interesting, yet they are quite complicated. Indeed, the observer-based construction includes 
the Atiyah-Sutcliffe constructions when setting $d=1$, and these have proved to be quite difficult problems to tackle (indeed, at the time or writing, very little is known when $n \geq 5$). We 
also remark that the two constructions for $X = \mathbb{C}\mathbb{P}^1$ are quite classical in character. Indeed, there are links between our constructions and the classical 
Lagrange interpolating polynomial, and generalizations of it to polynomials of several variables (and particularly with the Chung-Yao work in \cite{Chung-Yao}, since 
our constructions satisfy their geometric criterion). There are also links with the classical notion of the 
Veronese embedding. While the author did not find these constructions in the literature, yet he suspects that they might be known. In any case, the fact that linear independence holds for these two 
constructions in the case $X= \mathbb{C}\mathbb{P}^1$ gives weight to our conjectures, since the case $X= \mathbb{C}\mathbb{P}^1$ can be obtained from the hyperbolic versions of these two 
constructions in the limiting case as the configuration of points approaches $n$ distinct points on the sphere at infinity.

\def\Dbar{\leavevmode\lower.6ex\hbox to 0pt{\hskip-.23ex \accent"16\hss}D}
\providecommand{\bysame}{\leavevmode\hbox to3em{\hrulefill}\thinspace}
\providecommand{\MR}{\relax\ifhmode\unskip\space\fi MR }
\providecommand{\MRhref}[2]{%
  \href{http://www.ams.org/mathscinet-getitem?mr=#1}{#2}
}
\providecommand{\href}[2]{#2}

\end{document}